\documentclass{article}
\usepackage[normalem]{ulem}
\usepackage{datetime}
\usepackage{amsthm}
\newtheorem{theorem}{Theorem}[section]
\newtheorem{lemma}[theorem]{Lemma}
\usepackage{mathtools}
\usepackage{esint}
\usepackage{thmtools}
\usepackage{enumitem}

\usepackage{xcolor}

\declaretheoremstyle[headfont=\normalfont]{normalhead}
\usepackage{color}
\usepackage{graphicx}
\usepackage{morefloats}
\usepackage{hyperref}
\usepackage{ amssymb }
\usepackage{bbm}
\usepackage{mathrsfs}
\usepackage{textcomp}
\DeclareMathOperator*{\esssup}{ess\,sup}
\DeclareMathOperator*{\essinf}{ess\, inf}

\usepackage{url}
\usepackage{float}
\usepackage{biblatex}
\newtheoremstyle{mydef}
{\topsep}{\topsep}%
{}{}%
{\itshape}{}
{\newline}
{%
  \rule{\textwidth}{0.0pt}\\*%
  \thmname{#1}~\thmnumber{#2}\thmnote{\-\ #3}.\\*[-1.5ex]%
  \rule{\textwidth}{0.0pt}}%
  
\begin{document}

\theoremstyle{mydef}
\newcommand{\N}{\mathbb{N}}
\newcommand{\R}{\mathbb{R}}
\newcommand{\Hi}{\mathbb{H}}
\newcommand{\Q}{\mathbb{Q}}
\newcommand{\Z}{\mathbb{Z}}
\newcommand{\E} {\mathbb{E}}
\newcommand{\pv}{\textrm{P.V}}
\newcommand\numberthis{\addtocounter{equation}{1}\tag{\theequation}}
\newtheorem{definition}{Definition}
\newtheorem{conjecture}{Conjecture}
\newtheorem{question}{Question}
\newtheorem{result}{Result}
\newtheorem{remark}{Remark}
\newtheorem{proposal}{Proposal}
\newtheorem{corollary}{Corollary}
\newtheorem{observation}{Observation}

\date{}

\title{Helmholtz Solutions for the Fractional Laplacian and Other Related Operators}
\author{Vincent Guan \thanks{ Department of Mathematics, University of British Columbia, Vancouver, BC, V6T 1Z2, Canada; vguan23@math.ubc.ca},
Mathav Murugan\thanks{ Department of Mathematics, University of British Columbia, Vancouver, BC, V6T 1Z2, Canada; mathav@math.ubc.ca}, and Juncheng Wei\thanks{ Department of Mathematics, University of British Columbia, Vancouver, BC, V6T 1Z2, Canada; jcwei@math.ubc.ca} }

\maketitle
\begin{abstract}
\hskip -.2in
\noindent
We show that the bounded solutions to the fractional Helmholtz equation, $(-\Delta)^s u= u$ for $0<s<1$ in $\R^n$, are given by the bounded solutions to the classical Helmholtz equation $(-\Delta)u= u$ in $\R^n$ for $n \ge 2$ when $u$ is additionally assumed to be vanishing at $\infty$. When $n=1$, we show that the bounded fractional Helmholtz solutions are again given by the classical solutions $A\cos{x} + B\sin{x}$. We show that this classification of fractional Helmholtz solutions extends for $1<s \le 2$ and $s\in \N$ when $u \in C^\infty(\R^n)$. Finally, we prove that the classical solutions are the unique bounded solutions to the more general equation $\psi(-\Delta) u= \psi(1)u$ in $\R^n$, when $\psi$ is complete Bernstein and certain regularity conditions are imposed on the associated weight $a(t)$.

\end{abstract}
\section{Introduction}

\subsection{Motivation and Main Results}

\indent

An elementary result is that the bounded solutions to the equation $-u_{xx} = u$ in $\R$ are given in the form: $u(x)=A\cos{x} + B\sin{x}$. In higher dimensions, bounded solutions to the classical Helmholtz equation $-\Delta u = u$ on $\R^n$ are expressed in terms of Bessel functions and spherical harmonics (see Appendix and also \cite{Agmon}). These higher dimensional solutions are $C^\infty$, bounded, and vanishing at $\infty$. Helmholtz functions arise frequently when solving PDEs such as the heat equation and the wave equation via separation of variables. They also play fundamental roles in inverse problems and scattering theory. 

In recent years, there is a  growing interest in the fractional Laplacian $(-\Delta)^s$, with $0<s<1$, due to various applications  involving nonlocal diffusion \cite{2016}. A natural question is to classify  solutions to the fractional Helmholtz equation

\begin{equation}    (-\Delta)^s u(x) = u(x) \ \ \mbox{in} \ \R^n.
\tag{1.1}\label{eq: 1.1}
\end{equation}
In the first theorem of this paper, we prove that the bounded and vanishing fractional Helmholtz solutions to \eqref{eq: 1.1} are the same as the classical Helmholtz solutions.

\begin{theorem} 
\label{Th 1.1}
Let $u(x) \in L^\infty(\R^n)$. If $n=1$ and $0<s<1$, then all solutions to $(-\Delta)^s u = u$  in $\R$ are $ A \cos (x)+B \sin (x)$ for some constants $A,B \in \R$. If $n\geq 2$ and we assume that  $u(x) \to 0$ as $|x| \to \infty$, then $u$ solves $(-\Delta)^s u = u$  in $\R^n$ for $0<s<1$ if and only if $(-\Delta) u = u$ in $\R^n$.
\end{theorem}

\noindent
{\bf Remark:} The case of $n=1$ was first proved by Fall and Weth in \cite{FallWeth}.\\

We will then show that this result can be generalized to powers $1<s \le 2$ by reducing this case to the previous case $0<s<1$.

\begin{theorem}
\label{Th 1.2}
Let $u(x) \in C^\infty(\R^n) \bigcap L^\infty(\R^n)$. Then $u$ solves  $(-\Delta)^s u = u$  in $\R^n$ for $0<s \le 2$ if and only if $(-\Delta) u = u$ in $\R^n$. 
\end{theorem}

For $s>2$, we have the following classification result for the polyharmonic Helmholtz equation.

\begin{theorem}
\label{Th 1.2-2}
Let $u(x) \in C^\infty(\R^n) \bigcap L^\infty(\R^n)$ and $ m \in \N$.  Then $u$ solves  $(-\Delta)^m u = u$  in $\R^n$  if and only if $(-\Delta) u = u$ in $\R^n$. 
\end{theorem}

Finally, we will consider the Helmholtz equation for complete Bernstein functions $\psi$ of the Laplace operator. The Bernstein Helmholtz equation is given by

\begin{align*}
    \psi(-\Delta) u(x) = u(x) \ \ \mbox{in} \ \R^n.
    \tag{1.2}\label{eq: 1.2}
\end{align*}
$\psi(t)$ is said to be a complete Bernstein function if $\psi(t)=\frac{\mathcal{L}\{f\}(t)}{t}$, where $\mathcal{L}$ is the Laplace transform and $f: [0, \infty) \to [0, \infty)$ is completely monotone, i.e. $(-1)^k f^{(k)} \ge 0$. In particular, $\psi(t)$ generalizes $t^s$ with $0<s<1$.

\begin{theorem}
\label{Th 1.3}
Let $n\geq 2$. Suppose $u(x) \to 0$ as $|x| \to \infty$ and that the associated weight $a(t)$ for $\psi(-\Delta)$ in the extension problem \eqref{eq: 1.5} is $A_2$ and obeys $a(t) \sim t^{\alpha}$ for $t\gg 1$ and $-1<\alpha<1$. Then $u \in L^\infty(\R^n)$ solves $\psi(-\Delta) u = \psi(1)u$ in $\R^n$ if and only if $(-\Delta) u = u$ in $\R^n$.
\end{theorem}

\newpage

\subsection{Preliminaries}

\indent

Recall that $(-\Delta)^s$ is defined via the Fourier multiplier $|\xi|^{2s}$, such that in the distributional sense, we have

$$
\widehat{(-\Delta)^s u} (\xi) = |\xi|^{2s}\hat{u}(\xi).
$$
We will also use two other definitions. The first is via the singular integral $(-\Delta)^s u(x) = c_{n,s} \textrm{P.V} \int_{\R^n} \frac{u(x)-u(y)}{|x-y|^{n+2s}}dy$ where $c_{n,s}$ is a normalization constant. The second definition, made famous by Caffarelli and Silvestre in \cite{cs}, characterizes $(-\Delta)^s u(x)$ by its harmonic extension $\dot{u}(x,t)$. In particular, $\dot{u}(x,t) \in H^1_{loc}(\R^{n+1}_+, t^{1-2s})$ weakly solves the harmonic equation:

\begin{align*}
    & \nabla \cdot (t^{1-2s} \nabla u) = t^{1-2s}[ \Delta_x \dot{u} + \frac{1-2s}{t} \dot{u}_t + \dot{u}_{tt}] = 0 \textrm{ on } \R^{n+1}_{+}, \\
    & \lim_{t \to 0} t^{1-2s} \dot{u}_t(x,t) = -c_{n,s}(-\Delta)^s u \textrm{ on } \R^n, \\
    & \dot{u}(x,0) = u(x) \textrm{ on } \R^n \times \{0\}.
    \tag{1.3} \label{eq: 1.3}
\end{align*}

Similarly, for complete Bernstein functions $\psi(x)$, $\psi(-\Delta)$ is defined via the Fourier multiplier $\psi(|\xi|^2)$, such that in the distributional sense, we have

$$
\widehat{\psi(-\Delta) u}(\xi) = \psi(|\xi|^2)\hat{u}(\xi).
$$
From \cite{kwas}, for any complete Bernstein function $\psi$, we may characterize $\psi(-\Delta)$ using the harmonic extension $\dot{u}(x,t) \in H^1_{loc}(\R^{n+1}_+, a(t))$, which solves  
 \begin{align*}
    & \partial_s ^2 \dot{u}(x,s) + A(s) \Delta_x \dot{u}(x,s) = 0 \textrm{ on } \R^{n+1}_{+}, \\
    & \partial_s \dot{u}(0,x) = -\psi(-\Delta) u(x) \textrm{ on } \R^n,\\
    &\dot{u} (x,0)= u(x)  \textrm{ on } \R^n \times \{0\},
    \tag{1.4} \label{eq: 1.4} 
\end{align*}
where $A(s)$ is a locally finite measure on some interval $[0, R)$,
with $R$ possibly infinite. If $A(s)$ is non-negative and $L^1_{\textrm{loc}}$, then we may apply the change of variables: $ds = (a(t))^{-1}dt$ and $A(s)ds = a(t)dt$. We call $a(t)$ the associated weight to equation \eqref{eq: 1.2}. This reduces \eqref{eq: 1.4} to the familiar form:
\begin{align*}
    & \nabla \cdot (a(t) \nabla \dot{u}) =  0 \textrm{ on } \R^{n+1}_{+}, \\
    & \lim_{t \to 0} a(t) \dot{u}_t(x,t) = -c_{n,a}\psi(-\Delta) u(x) \textrm{ on } \R^n. 
    \tag{1.5} \label{eq: 1.5}
\end{align*}
 For the rest of the paper, the trace condition $\dot{u} (x,0)= u(x)$ will be implicitly assumed.

\subsection{Comments about the proofs}

\indent

Since  $(-\Delta)^s$ and $\psi(-\Delta)$ are non-local operators,  it is difficult to solve these nonlocal Helmholtz equations directly. We cannot apply separation of variables as in the $s=1$ case (shown in Appendix), and in general, it is not easy to compute the integral $(-\Delta)^su$ explicitly. To prove compatibility of classical solutions with other Helmholtz equations, it is helpful to consider the equivalent extension problem, i.e. \eqref{eq: 1.3} for $(-\Delta)^s$ and \eqref{eq: 1.4}, \eqref{eq: 1.5} for $\psi(-\Delta)$. This way, the problem is reformulated as a second order elliptical partial differential equation.

Uniqueness of fractional Helmholtz solutions in the class of bounded vanishing functions for $n \ge 2$ is achieved by working with the extension problem and using standard energy techniques, following the seminal work of \cite{Frank}.  The $n=1$ case is handled separately.


Similar but more sophisticated techniques will be used to prove uniqueness in the general case where we consider the solutions to $\psi(-\Delta)u = u$. $a(t) \in A_2$ is assumed to achieve appropriate H\"{o}lder estimates and the asymptotic behaviour $a(t) \sim t^\alpha$ is assumed to apply the estimate $|\dot{u}_t| \le \frac{C}{t}$ from \cite[Prop 4.6]{cabre}.

\section{Fractional Helmholtz Solutions $0<s<1$}

\textbf{Proof of Theorem \ref{Th 1.1}}: When $n \ge 2$, Theorem \ref{Th 1.1} is a corollary of Lemma \ref{lm 2.1} and Lemma \ref{lm 5.2}. When $n=1$, Theorem \ref{Th 1.1} is a corollary of Lemma \ref{lm 2.1} and Lemma \ref{lm 2.3}. \qed

\begin{lemma}
\label{lm 2.1}
If $u\in L^\infty(\R^n)$ solves $-\Delta u = u $ in $\R^n$, then  $(-\Delta)^s u = u$ in $\R^n$ for $0<s<1$.
\end{lemma}

\noindent
\textbf{Proof}: Suppose that $u\in L^\infty(\R^n)$ and $-\Delta u = u$ on $\R^n$. It suffices to show that the extension problem \eqref{eq: 1.3} for $(-\Delta)^s u = u$ is solvable. Hence, it suffices to show that there is an extension  $\dot{u}(x,t)=u(x)\phi(t)$ with $\phi(0)=1$ such that
\begin{align*}
    & \phi(t)\Delta u(x)+ u(x)\frac{1-2s}{t}\phi'(t)+u(x)\phi''(t)=0 \textrm{ on } \R^{n+1}_+, \\
    & \lim\limits_{t \to 0} t^{1-2s}u(x)\phi'(t) = -c_{n,s}(-\Delta)^s u(x) = -c_{n,s}u(x) \textrm{ on } \R^n. 
    \tag{2.1} \label{eq: 2.1}
\end{align*}
Since $u$ is a classical Helmholtz solution, \eqref{eq: 2.1} reduces to
\begin{align*}
    -\phi(t) + \frac{1-2s}{t}\phi'(t)+\phi''(t)=0 \textrm{ for } t \ge 0, \\
    \lim\limits_{t \to 0} t^{1-2s}\phi'(t) = -c_{n,s} \in \mathbb{R}.
    \tag{2.2} \label{eq: 2.2}
\end{align*}
It suffices to solve \eqref{eq: 2.2} with the initial condition $\phi(0)=1$. This ODE is solved in  \cite[Section 3.2]{cs}. 
\qedsymbol \\
    
We now show that \cite[Thm 1]{Frank}  can be generalized to non-radial functions in the following sense:

\begin{lemma}
\label{lm 2.2}
Let $n \ge 2$, $V(r) \in C^\infty(\R^n)$ be radially non-decreasing and $u$ be a bounded and vanishing function satisfying $(-\Delta)^s u+ V(r) u=0$ on $\R^n$. Then the projection of $u$ onto any spherical harmonical eigenfunction is unique up to a constant factor.
\end{lemma}

\noindent
\textbf{Proof:}  Let $n \ge 2$ and $u \in L^\infty(\R^n)$ be vanishing and solve $(-\Delta)^s u + V(r)u= 0$ on $\R^n$ where $V(r)$ is radially non-decreasing and differentiable. First, we use eigenfunction decomposition with spherical harmonics to reduce this problem to the case where $u$ is additionally assumed to be radial. In particular, we consider the extension $\dot{u}(x, t)$, which solves \eqref{eq: 1.3}. Using spherical harmonics, we obtain the decomposition
\begin{align*}
    &\dot{u}(r, \theta, t) = \sum_l \dot{u_l}(r,t)\phi_{l}(\theta),
\end{align*}
where $\Delta_{\mathbb{S}^{n-1}}\phi_{l}(\theta) = 
-\mu_{l}\phi_{l}(\theta)$ and $\mu_l= l(l+n-2)$. Note that for each fixed eigenvalue $\mu_l$, there may be multiple associated eigenfunctions $\phi_l$ included in the sum. Then, since $\Delta_x u = \partial_r^2 u + \frac{n-1}{r} \partial_r u + \frac{ \Delta_{\mathbb{S}^{n-1}} u}{r^2}$ and $\{\phi_l\}_l$ is an orthogonal family, for each $l$, we may substitute $\dot{u}=\dot{u_l}(r,t)\phi_l(\theta)$ into \eqref{eq: 1.3} to yield

\begin{align*}
    & \dot{u}_{l, rr} + \frac{n-1}{r}\dot{u}_{l,r} + \frac{1-2s}{t}\dot{u}_{l,t} + \dot{u}_{l, tt} - \frac{\mu_l \dot{u}_l}{r^2} = 0 \textrm{ on } \R^{n+1}_+,\\
    & \lim\limits_{t \to 0} t^{1-2s}\dot{u}_{l,t}  = c_{n,s} V(r)\dot{u_l}(r,0) := c_{n,s}V(r)u_l(r) \textrm{ on } \R^n.
     \tag{2.3} \label{eq: 2.3}
\end{align*}
Note that $u(x) =\sum \dot{u}_l(r,0)\phi_{l}(\theta) =\sum u_l(r)\phi_{l}(\theta)$ by substituting $t=0$. Now, let $\dot{v}_l(r,t) = r^{-l} \dot{u}_l(r,t)$ and $v_l(r) = r^{-l} u_l(r)$. This yields the system:

\begin{align*}
    &\dot{v}_{l, rr} + \frac{2l+n-1}{r}{\dot{v}}_{l,r} + \frac{1-2s}{t}\dot{v}_{l,t} + {\dot{v}}_{l, tt}  = 0 \textrm{ on } \R^{n+1}_+,\\
    &\lim\limits_{t \to 0} t^{1-2s}\dot{v}_{l,t}  = c_{n,s}V(r){v}_l(r) \textrm{ on } \R^n.
\tag{2.4} \label{eq: 2.4}
\end{align*}
In particular, the last term of the first condition of \eqref{eq: 2.3} cancels since $\mu_l = l(l+n-2) = l(l-1)+l(n-1)$. Now, note that each $v_l$ is bounded, radial, and vanishing, so we may use the proof of \cite[Thm 1]{Frank} to show that each $v_l$ is unique. To simplify notation, we redefine $u=v_l$, $\dot{u}=\dot{v_l}$, and $V(r)=c_{n,s}V(r)$. Following \cite[(4.6)]{Frank}, we define the energy $H(r)$ on $u$ as
\begin{equation}
    H(r) = \frac{1}{2}\left[\int_{0}^{\infty}t^{1-2s}(\dot{u}_{r}^2-\dot{u}_{t}^2)dt -V(r) u (r)^2\right].
    \tag{2.5} \label{eq: 2.5}
    \end{equation}
In \cite[Prop B.2]{Frank}, it is shown that $H(0) \le -\frac{1}{2} V(0)u(0)^2$ and $H(\infty) = 0$  for any radial $u(r)$ vanishing at infinity. We now show that $H'(r) \le 0$ for $u(r)=v_l(r)$. We compute
\begin{align}
H'(r) &= \int_{0}^{\infty}\frac{d}{dr}\frac{t^{1-2s}}{2}(\dot{u}_r^2-\dot{u}_t^2)dt  - \frac{1}{2}V'(r)u(r)^2 - V(r)u(r)u'(r) \nonumber \\
                &= \int_{0}^{\infty}t^{1-2s}\dot{u}_r \dot{u}_{rr} dt - \int_{0}^{\infty}t^{1-2s}\dot{u}_t \dot{u}_{rt}dt -  \frac{1}{2}V'(r)u(r)^2 - V(r)u(r)u'(r) ) \nonumber \\ 
                &= \int_{0}^{\infty}t^{1-2s}\dot{u}_r \dot{u}_{rr}dt - t^{1-2s}\dot{u_t}\dot{u_r} \vert_{t=0}^{t=\infty}  \nonumber   \\
                 & \qquad   + \int_{0}^{\infty}t^{1-2s}\dot{u}_r \left(\frac{1-2s}{t}\dot{u}_{t} + \dot{u}_{tt}\right)dt - \frac{1}{2}V'(r)u(r)^2 - V(r)u(r)u'(r) ) \nonumber \\
                & =  \int_{0}^{\infty}t^{1-2s}\dot{u}_r\left(\dot{u}_{rr} + \frac{1-2s}{t}\dot{u}_t + \dot{u}_{tt}\right)dt  - \frac{1}{2}V'(r)u(r)^2 ) \nonumber \\
                & = -\frac{2l+n-1}{r}\int_{0}^{\infty}t^{1-2s}\dot{u}_r^2dt-\frac{1}{2}V'(r)u(r)^2 \le 0.
\tag{2.6} \label{eq: 2.6}
\end{align}
In the second to last equality, note that $-t^{1-2s}\dot{u_t}\dot{u_r} \vert_{t=0}^{t=\infty} = V(r)u(r)u'(r)$ results from applying the boundary condition \eqref{eq: 2.4} and the estimate $|\dot{u_t}(x,t)| \le \frac{C}{t}$ from \cite[Prop 4.6]{cabre}. The last equality is due to the first condition in \eqref{eq: 2.4} and this expression is non-positive since $n \ge 2, l \ge 0$, and $V'(r) \ge 0$.

If $u(0)=0$, then it is clear that $H(r)=0$ and $\frac{dH}{dr}=0$ given the properties of $H(r)$ established above. If $n \ge 2$, then we may immediately conclude that $u_r=0$ from \eqref{eq: 2.6} and hence $u=0$. Equivalently, since \eqref{eq: 2.4} is a linear system, $u=v_l$ must be the unique bounded and vanishing solution to \eqref{eq: 2.4} up to a constant factor.
\qed

\begin{lemma}
\label{lm 2.3}
    If $u \in L^\infty(\R)$ solves $(-\Delta)^s u = u $ in $\R$, then $u^{''} = u $ in $\R$.
\end{lemma}

\noindent
\textbf{Proof}: Recall that the Helmholtz solutions on $\R$ are $u(x) = A\cos(x) + B\sin(x)$. We may assume $(-\Delta)^s u = u$ and $u \in S^{'}(\R)$, the space of tempered distributions. Then, taking the Fourier transform on both sides gives $|\xi|^{2s}\hat{u}(\xi) = \hat{u}(\xi)$ in the weak sense. A quick proof is shown below. For all $\phi \in S(\R)$:

\begin{equation}
    \langle \hat{u}, \phi \rangle = \int{((-\Delta)^s u) \hat{\phi}} = \int{u (-\Delta)^s \hat{\phi}} = \int{\hat{u} |\xi|^{2s}\phi} =\langle |\xi|^{2s}\hat{u}(\xi), \phi \rangle
    \tag{2.7} \label{eq: 2.7}
\end{equation}

Because the Fourier transform is a bijective map from $S^{'} \to S^{'}$, we note that $\hat{u}(\xi) \in S^{'}(\R)$ and it has support $\{-1, 1\}$. Thus, we may write $\hat{u}(\xi)= \sum_{\{\alpha \le N\}}C_\alpha \partial^\alpha \delta_{-1} + D_\alpha \partial^\alpha \delta_{1}$ for some constants $N, C_{\alpha}, D_\alpha$. Then, taking the inverse transform, we see that $u(x) = a(x)\sin(x) + b(x)\cos(x)$ for some polynomials $a(x)$ and $b(x)$.  Since $u$ is bounded,  $u(x)=A\cos(x)+B\sin(x)$. \qed

\section{Fractional Helmholtz Solutions $1<s\le2$}

\indent

In this section, we consider fractional Helmholtz solutions with $ 1<s\le2$. Existence can be proved in a similar way as in Lemma \ref{lm 2.1}, by using  the extension problem from \cite{Yang}. Instead we present a new proof by decomposition of nonlocal operators.  To do so, we use the fact that we can rewrite $(-\Delta)^s u$ with $1<s\le2$ in terms of the standard fractional laplacian with $\frac{1}{2}<s/2\le1$ when $u$ has sufficient regularity. We also use this decomposition to prove uniqueness.\\

\noindent
\textbf{Proof of Theorem \ref{Th 1.2}}: This is a corollary of Lemma \ref{lm 3.2} and Lemma \ref{lm 3.3}. \qed

\begin{lemma}
\label{lm 3.1}
Let $u \in C^\infty(\R^n) \bigcap L^\infty(\R^n)$. Then for $ 1<s\le2$, we have $(-\Delta)^s u = (-\Delta)^{s/2} [(-\Delta)^{s/2}u]$.
\end{lemma}

\noindent
\textbf{Proof:} Let $u \in C^\infty(\R^n) \bigcap L^\infty(\R^n)$ and $1<s\le2$. We consider the definition of $(-\Delta)^s$ from \cite{powers}:
$$ (-\Delta)^s u = L_{2,s}u  =  c_{n,2,s}\int \frac{u(x-2y)-4u(x-y)+6u(x)-4u(x+y)+u(x+2y)}{|y|^{n+2s}}dy.$$
Note also that $L_{1,s/2}u = (-\Delta)^{s/2}u$. Thus, we want to show that  $L_{1,s/2}(L_{1,s/2}u) = L_{2, s}u$. In \cite[Thm 1.9]{powers} this is obtained for all $\phi \in C_c^\infty(\R^n)$ since
\begin{equation}
    \mathscr{F}(L_{2,s}\phi)(\xi) =|\xi|^{2s} \mathscr{F}(\phi)(\xi)= |\xi|^s  \mathscr{F}(L_{1,s/2}\phi )(\xi) = \mathscr{F}(L_{1,s/2}(L_{1,s/2}\phi))(\xi).
\tag{3.1} \label{eq: 3.1}
\end{equation} 

Now, since $u \in L^\infty(\R^n) \bigcap C^\infty(\R^n)$, we may apply \cite[Lemma 1.5]{powers}, which tells us that for any $\phi \in C_c^\infty(\R^n)$,

\begin{align*}
    & \int_{\R^n}u(x)L_{2,s}\phi(x)dx =  \int_{\R^n}\phi(x)L_{2,s}u(x)dx, \\
    & \int_{\R^n}u(x)L_{1,s/2}\phi(x)dx =  \int_{\R^n}\phi(x)L_{1,s/2}u(x)dx.
\tag{3.2} \label{eq: 3.2}
\end{align*}

The desired result then follows from applying \eqref{eq: 3.1} and \eqref{eq: 3.2}.

\begin{align*}
    & \int_{\R^n}\phi(x)L_{2,s}u(x)dx = \int_{\R^n}u(x)L_{2,s}\phi(x)dx = \int_{\R^n}u(x)L_{1,s/2}(L_{1,s/2}\phi(x))dx \\
    &= \int_{\R^n}L_{1,s/2}
    u(x)L_{1,s/2}\phi(x)dx = \int_{\R^n} \phi(x)L_{1,s/2}(L_{1,s/2}u(x))dx
\tag{3.3} \label{eq: 3.3}
\end{align*}
\qed



\begin{lemma}
\label{lm 3.2}
If $u$ solves $(-\Delta) u = u $ in $\R^n$ and $u \in C^\infty(\R^n) \bigcap L^\infty(\R^n)$, then $u$ also solves  $(-\Delta)^s u = u$ in $\R^n$ for $1<s\le2$.
\end{lemma}

\noindent
\textbf{Proof:}
Let $1<s \le 2$ and $u \in C^\infty(\R^n) \bigcap L^\infty(\R^n)$ solve $(-\Delta) u = u $ in $\R^n$. By Theorem \ref{Th 1.1}, $(-\Delta)^{s/2} u = u$ since $s/2 \in (0,1]$. We then apply Lemma \ref{lm 3.1} to yield

\begin{equation}
    (-\Delta)^s u = (-\Delta)^{s/2}[(-\Delta)^{s/2}u] = (-\Delta)^{s/2}u = u.
    \tag{3.4} \label{eq: 3.4}
\end{equation} \qedsymbol

To prove uniqueness of these fractional solutions, it suffices to show that all solutions $(-\Delta)^s u = u$, with $1<s \le 2$, satisfy $(-\Delta)^{s/2} u = u$, since uniqueness has already been shown for $s \in (0,1]$ with Lemma \ref{lm 2.2} and Lemma \ref{lm 2.3}.

\begin{lemma}
\label{lm 3.3}
Let $1<s\le2$. If $u \in C^\infty(\R^n) \bigcap L^\infty(\R^n)$ solves $(-\Delta)^s u = u$ on $\R^n$, then $(-\Delta)^{s/2} u = u$ on $\R^n$.
\end{lemma}

\noindent
\textbf{Proof:}
Let $1<s\le2$. and $(-\Delta)^s u  = u$ on $\R^n$. We can define the system
\begin{equation}
    \left\{\begin{array}{l}
v=(-\Delta)^{s/2}u-u, \\
(-\Delta)^{s/2}v+v=0.
\end{array}
    \right.
      \tag{3.5} \label{eq: 3.5}
\end{equation}

By Lemma \ref{lm 3.1}, it suffices to show that $v=0$. 
We proceed by using the Maximum Principle. Let $w(x) = C+ |\sqrt{1+x^2}|^\rho$ with $\rho < s$. We first claim that we can construct $w$ such that $(-\Delta)^{s/2}w+w>0$. To see this, note that $ |\sqrt{1+x^2}|^\rho$ asymptotically behaves like $|x|^\rho$  while remaining uniformly bounded for $x$ close to $y$. Hence, we may let $x = \alpha \tilde{x}, y= \alpha \tilde{y}$ and compute

\begin{align*}
& (-\Delta)^{s/2}w(x) = c_{n,s} \pv  \int_{\R^n} \frac{ |\sqrt{1+x^2}|^\rho- |\sqrt{1+y^2}|^\rho}{|x-y|^{n+s}}dy \\
&= c_{n,s} \pv \int_{\R^n} \frac{ |\sqrt{1+\alpha^2\tilde{x}^2}|^\rho- |\sqrt{1+\alpha^2 \tilde{y}^2}|^\rho}{\alpha^{n+s}|\tilde{x}-\tilde{y}|^{n+s}} \alpha^n d\tilde{y} \sim \alpha^{\rho-s} (-\Delta)^s w(\tilde{x}).
\tag{3.6} \label{eq: 3.6}
\end{align*}
Then, since $\rho-s<0$ and $(-\Delta)^s w(\tilde{x})$ is bounded because the integrand in the last line of \eqref{eq: 3.6} is $o(|y|^{n})$, we may pick $\alpha, \rho, C$ such that $(-\Delta)^{s/2}w+w>0$. Now, let $\Tilde{v} = v- \epsilon w$. We immediately see that

\begin{align*}
    &(-\Delta)^{s/2}\Tilde{v} + \Tilde{v} = (-\Delta)^{s/2} v + v - \epsilon((-\Delta)^{s/2}w+w) \\
    &= - \epsilon((-\Delta)^{s/2}w+w) < 0.
    \tag{3.7} \label{eq: 3.7}
\end{align*}
Note that $v = o(w)=o(|x|^\rho)$ since we can choose $\rho$ arbitrarily close to $s$ and $v = o(|x|^s)$. Therefore, $\Tilde{v}(\infty) = -\infty$ and $\tilde{v} \in C^2(\R^n)$ by Lemma \ref{lm 5.4}. From this, we may deduce that there is some $x_0$ such that $\tilde{v}(x_0)=\max\limits_{x \in \R^n}v(x)$. Hence, we may use the fact that $(-\Delta)^{s/2} \Tilde{v}(x_0) = \textrm{P.V} \int \frac{\Tilde{v}(x_0) - \Tilde{v}(y)}{|x_0-y|^{n+s}} \ge 0$ and  \eqref{eq: 3.7} to obtain the inequalities

\begin{equation}
    0 \le (-\Delta)^{s/2}\Tilde{v}(x_0) < -\Tilde{v}(x_0) \implies v(x) - \epsilon w(x) \le v(x_0) - \epsilon w(x_0) \le 0.
    \tag{3.8} \label{eq: 3.8}
\end{equation} 
Taking $\epsilon \to 0$ yields $v(x) \le 0$. By a symmetric argument, if we let $\Tilde{v} = \epsilon w - v$, then we would have $(-\Delta)^{s/2}\Tilde{v} + \Tilde{v} > 0$ and $v(x) \ge 0$. This proves $v(x)= 0$ as desired. \qed

\section{Polyharmonic Helmholtz Solutions}

\indent

In this section, we consider Helmholtz solutions to the equation 

\begin{equation}
    (-\Delta)^m u = u \textrm{ in } \R^n, \textrm{ where } m \in \N. 
     \tag{4.1} \label{eq: 4.1}
\end{equation}

\noindent
\textbf{Proof of Theorem \ref{Th 1.3}}: From Lemma \ref{lm 3.2}, it is clear that $-\Delta u = u$ will solve $(-\Delta)^m u = u$ for any $m \in \N$ and we have already shown uniqueness when $m=1,2$. Let us now consider the case $m=3$. To prove uniqueness of the classical solutions in the class of bounded functions, we may consider the system

\begin{align*}
    v= -\Delta u - u,\\
    \Delta^2 v - \Delta v + v=0.
    \tag{4.2} \label{eq: 4.2}
\end{align*}

It suffices to show that $v=0$. To do so, let $\eta$ be a cutoff function supported in $B_{2R}$. We may multiply the bottom equation by $\eta^2 v$ and integrate by parts to yield

\begin{align*}
    &\int (\Delta (\eta v))^2 + \int |\nabla \eta v|^2 + \int \eta^2 v^2 = -2 \int v \Delta v |\nabla \eta|^2 + v\eta \nabla v \cdot \nabla \eta  \\ 
    & + 4\int (\nabla \eta \cdot \nabla v)^2 + v\Delta \eta \nabla \eta \cdot \nabla v + \int v^2(\Delta \eta)^2 - v^2 \eta \Delta \eta.
    \tag{4.3} \label{eq: 4.3}
\end{align*}

Then, note that we may pick $\eta$ such that $|\nabla \eta| \lesssim \frac{1}{R}$. Also, from standard elliptical theory, since $u \in L^\infty(\R^n)$, it follows that all of its derivatives are uniformly bounded as well. Thus, after using integration by parts and Cauchy-Schwarz on the right hand side of \eqref{eq: 4.3}, we obtain the inequality

\begin{align*}
    &\int_{B_{2R}} (\Delta v)^2 +|\nabla v|^2 + v^2  \lesssim \frac{1}{R^2} \int_{B_{2R}}v^2 + |\nabla v|^2.
    \tag{4.4}\label{eq: 4.4}
\end{align*}

Hence, we may iterate the inequality \eqref{eq: 4.4} on its own right hand side, such that the bound is scaled by $\frac{1}{R^2}$ each time. Taking $R \to \infty$ then yields the desired result $v \equiv 0$.

When $m\geq 4$ we modify the system \eqref{eq: 4.2} such that $ v= -\Delta u- u$ satisfies $ \sum_{j=0}^{m-1} (-\Delta)^j v + v=0$. The case of $m\geq 4$ can be proved similarly as that of $m=3$. We omit the details. This proves Theorem \ref{Th 1.2-2}.

\qed

\section{Complete Bernstein Helmholtz Solutions}

We now classify solutions to the Helmholtz equation $\psi(-\Delta)u = u$ in $\R^n$ when $\psi$ is complete Bernstein. Without loss of generality, we assume that the extension $\dot{u}$ is bounded in $\R^{n+1}_+$.\\

\noindent
\textbf{Proof of Theorem  \ref{Th 1.3}}: Theorem \ref{Th 1.3} is a corollary of Lemma \ref{lm 5.1} and Lemma \ref{lm 5.2}. \qed

\begin{lemma}
\label{lm 5.1}
If $-\Delta u = u$ on $\R^n$, then $u$ also solves $\psi(-\Delta) u = \psi(1) u$ on $\R^n$.
\end{lemma}


\noindent
\textbf{Proof:} Recall the extension problem for $\psi(-\Delta) u = \psi(1)u$ on $\R^n$ given by 
\begin{align*}
    &\partial_s ^2 \dot{u}(x,s) + A(s) \Delta_x \dot{u}(x,s) = 0 \textrm{ on } \R^{n+1}_+, \\
    &\partial_s \dot{u}(0,x) =  -\psi(1) u(x) \textrm{ on } \R^n.
     \tag{5.1} \label{eq: 5.1}
\end{align*}
Now, let $\dot{u}(x,s) = \phi(s) u(x)$ such that $\phi(0)=1$. Since $-\Delta u(x) = u(x)$ on $\R^n$, the problem reduces to finding a function $\phi(s)$ satisfying

\begin{align*}
    &\phi''(s) = A(s) \phi(s) \textrm{ for } s \ge 0, \\
    &\phi'(0)=-\psi(1) \in \R,\\
    &\phi(0)=1.
     \tag{5.2} \label{eq: 5.2}
\end{align*}
By \cite[Section 3.1]{kwas}, a solution $\phi$ exists. In particular, we may let $\lambda=1$ in \cite[(3.1)]{kwas} and set $\varphi_1=\phi $. \qed \\


To show uniqueness, we have only been able to provide a proof in the special case where $a(t) \in A_2$ and $a(t) \sim t^{\alpha}$, $|\alpha|<1$ for $t\gg 1$. By scaling  the solution, it suffices to consider the case $\psi(1)=1$.

\begin{lemma}
\label{lm 5.2}
Let $u$ be bounded and vanishing solution to the  equation $\psi(-\Delta) u = \psi (1) u$ on $\R^n$ where $n \ge 2$ and the associated weight $a(t)$ in the extension problem \eqref{eq: 1.5} is $A_2$ and obeys $a(t) \sim t^{\alpha}$, $|\alpha|<1$ for $t\gg 1$. Then $u$ satisfies $ -\Delta u=  u $.
\end{lemma}

\noindent
\textbf{Proof:} Let $\psi(-\Delta)u = \psi (1) u$ on $\R^n$. We mimic \eqref{eq: 2.2}-\eqref{eq: 2.3} from the proof of Lemma \ref{lm 2.2}. Consider the extension $\dot{u}(x,t)$, which solves \eqref{eq: 1.5}, and decompose it with spherical harmonics to express $\dot{u}(r, \theta,t) = \sum \dot{u}_l(r,t)\phi_{l}(\theta)$. Then, let $\dot{v}_l(r,t) = r^{-l} \dot{u}_l(r,t)$ and $v_l = r^{-l} u_l(r) := r^{-l} \dot{u}_l(r,0)$. For each $l$, $\dot{v}_l$ satisfies
\begin{align*}
     &\partial_t(a(t) \partial_r \dot{v}_l) + \frac{2l+n-1}{r}a(t)\partial_r \dot{v}_l+a(t) \partial_r^2 \dot{v}_l = 0,\\
     & \lim_{t \to 0} a(t) \partial_t  \dot{v}_l(r,t) = -c_{a,s}\psi (1) v_l(r).
\tag{5.3} \label{eq: 5.3}
\end{align*}

To prove that each $v_l$ is unique up to a constant factor in the class of bounded and vanishing functions, it suffices to show that $v_l(0)=0$ implies  $v_l \equiv 0$, since \eqref{eq: 5.3} is a linear system. For simplicity, let $\dot{v}_l = \dot{v}$, $v_l = v$, and $-c_{a,s}\psi (1)= c$. Then, multiply the first equation in (\ref{eq: 4.3})  by $\dot{v}_r$ to get
\begin{align*}
    &\partial_t(a(t) \partial_r \dot{v})\partial_r \dot{v} + \frac{2l+n-1}{r}a(t)(\partial_r \dot{v})^2+a(t) \partial_r \dot{v}\partial_r^2 \dot{v}=0.
        \tag{5.4} \label{eq: 5.4}
\end{align*}
Note that the middle term is non-negative. We now integrate the other two terms over $(0, \infty) \times (0, \infty)$ and show that these integrals are also non-negative.
\begin{align*}
 &\int_{0}^{\infty}\int_{0}^{\infty}(a(t) \dot{v}_t)_t  \dot{v}_r + a(t) \dot{v}_r  \dot{v}_{rr} dtdr\\
&= \int_{0}^{\infty}\int_{0}^{\infty} a(t) \left(\frac{v_r^2-v_t^2}{2} \right)_r drdt + \int_{0}^\infty a(t)\dot{v}_t\dot{v}_r\vert_{t=0}^{t=\infty} dr \\
&= \int_{0}^{\infty} \frac{a(t)}{2} \left(\dot{v}_r^2(r,t)- \dot{v}_t^2(r,t)\right)\vert_{r=0}^{r=\infty}dt + \int_{0}^\infty \lim_{t \to \infty} a(t)\dot{v}_t(r,t) \dot{v}_r(r,t) + cv(r)\dot{v}_r(r,0) dr \\
&= \int_{0}^{\infty} \frac{a(t)}{2} \left( \lim_{r \to \infty} \dot{v}_r^2(r,t) + \dot{v}_t^2(0,t)\right)dt + c\int_{0}^\infty \left(\frac{\dot{v}^2(r,0)}{2}\right)_r dr \\
&= \int_{0}^{\infty} \frac{a(t)}{2} \left( \lim_{r \to \infty} \dot{v}_r^2(r,t) + \dot{v}_t^2(0,t)\right)dt - \frac{cv^2(0)}{2}.
 \tag{5.5} \label{eq: 5.5}
\end{align*}
Three cancellations were performed to obtain the second to last equality. First, $\lim\limits_{r \to \infty}a(t)\dot{v}_t(r, t)^2=0$ due to Parseval's Theorem and Lemma \ref{lm 5.3}. Second, $a(t)\dot{v}_r(0, t)^2=0$ since $\dot{v}_r(0,t)=0$ follows from Lemma \ref{lm 5.4}. Third, note that $\lim_{t \to \infty} a(t)\dot{v}_t(r,t) \dot{v}_r(r,t)=0$ follows from the estimate $|\dot{v}_t(x,t)| \le \frac{C}{t}$, which can be shown from a  similar rescaling argument done in the proof of \cite[Prop 4.6]{cabre}. This is where we use the assumption that for sufficiently large $t$, $a(t) \sim t^\alpha$ for some $|\alpha|<1$. Then, from \eqref{eq: 5.4} and \eqref{eq: 5.5}, we may easily compute


\begin{align*}
    \frac{1}{2}\int_{0}^{\infty} a(t) (\dot{v}_r(\infty, t)^2 + \dot{v}_t(0, t)^2) dt+ \int_{0}^{\infty}\int_{0}^{\infty} \frac{2l+n-1}{r}a(t) (\dot{v}_r)^2 drdt  - \frac{cv^2(0)}{2}= 0.
    \tag{5.6} \label{eq: 5.6}
\end{align*}
Each of these terms are non-negative since $a(t) \in A_2$ implies 
that $a(t)> 0$ a.e. and $v(0)=0$ is assumed. Thus, if $n \ge 2$, then we may conclude $\dot{v}_r(r,t) = 0$ $\forall t \ge 0$ . Since $v(0)=0$, then $v \equiv 0$ as desired. Hence, each $u_l(r)$ is unique up to a constant factor and $u(r, \theta)= \sum_l c_l u_l(r)\phi_l(r)$ for some constants $c_l$. The result follow from Lemma \ref{lm 5.1} and the fact that if $u(r, \theta) = \sum_l u_l(r)\phi_l(r) $ solves $-\Delta u = u$, then $\tilde{u}(r, \theta)= \sum_l c_l u_l(r)\phi_l(r)$ also solves $-\Delta u = u$.

\qed


\begin{lemma}
\label{lm 5.3}
Let $u \in L^\infty(\R^n)$ solve \eqref{eq: 1.2} and $\dot{u}$ solve \eqref{eq: 1.5} with $a(t) \in A_2$ and $\dot{u}(x,0)$ vanishing at infinity. Then we have $\lim\limits_{|x| \to \infty}a(t)\dot{u}_t(x, t)^2=0$.
\end{lemma}

\noindent
\textbf{Proof:}
Let $\Omega=B_R(x_0) \times [0,R] \subset \R^{n+1}_+$. By expanding and rearranging the first condition of the extension problem \eqref{eq: 1.5}, $\dot{u}$ solves 

\begin{align*}
    (a(t) \dot{u}_t)_t= -a(t)\Delta_x \dot{u} \textrm{ in } \Omega, \\
    a(t)\dot{u}_t=u \textrm{ on }  \Omega \cap \{t=0\}.
    \tag{5.7}\label{eq: 5.7}
\end{align*}
This holds for any fixed $x_0$ and $R$. Let $(x,t) \in \Omega$. By the Fundamental Theorem of Calculus,
\begin{align*}
    a(t)\dot{u}_t(x,t) - u(x) = - \int_{0}^{t}a(s)\Delta_x \dot{u}(x,s)ds.
    \tag{5.8}\label{eq: 5.8}
\end{align*}
Since $a(t)$ is locally integrable, we obtain the estimate 
\begin{align*}
    & \|a(t)\dot{u}_t(x,t)\|_{L^\infty(\Omega)} \le  \|u\|_{L^\infty(\Omega)} + \left|\int_{0}^{R}a(t)dt\right| \|\Delta_x \dot{u}(x,t)\|_{L^\infty(\Omega)} \\
    & \le \|u\|_{L^\infty(\Omega)} + C\|\Delta_x \dot{u}(x,t)\|_{L^\infty(\Omega)},
    \tag{5.9}\label{eq: 5.9}
\end{align*}
where $\|u\|_{L^\infty(\Omega)} \to 0$ as $x_0 \to \infty$ since $u(x)$ is vanishing. The result then follows after applying interpolation inequalities \cite[(6.83), (6.85)]{gilbarg} with Lemma \ref{lm 5.4}. \qedsymbol


\begin{lemma}
\label{lm 5.4}
If $u \in L^\infty(\R^n)$ solves $\psi(-\Delta) u = u$ on $\R^n$ and $a(t)$, from the extension problem \eqref{eq: 1.5}, is an $A_2$ weight, then $\dot{u}(x, \cdot) \in C^\infty(\R^n)$ and hence $u \in C^\infty(\R^n)$.
\end{lemma}

\noindent
\textbf{Proof:} Let $Q_R = B_R(0) \times (0,R) \subset \R^{n+1}_+$ and $\partial' Q_R =  B_R(0) \times \{0\}$. The result follows after adapting \cite[Thm 1.2]{tan} for $a(t) \in A_2$ from the case $a(t)=t^{1-2s}$. That is, we want to show that if $\dot{u}$ solves \eqref{eq: 5.7} with $\Omega = Q_1$, then

\begin{align*}
\sup_{Q_{1/2}} \dot{u} \le C(\inf_{Q_{1/2}}\dot{u}).
\tag{5.10}\label{eq: 5.10}
\end{align*}

This is because \eqref{eq: 5.10} implies that $\|\dot{u}\|_{C^\beta(Q_1)}\le C\| \dot{u}\|_{L^\infty(Q_1)} \le K < \infty$ for some $\beta>0$ \cite[Thm 2.3.15]{fabes}. Below, in line \eqref{eq: 5.11}, we show that when $\beta<1$, $\dot{u}(x, \cdot) \in C^\beta(B_1)$ implies $\dot{u}(x, \cdot) \in C^{2\beta}(B_1)$ with a uniform bound independent of $t$. Hence, we know that $\dot{u}(x, \cdot) \in C^1(B_1)$ and we may then invoke \eqref{eq: 5.10} on $\nabla_x \dot{u}$ and repeat the doubling argument to show that $\dot{u}(x, \cdot) \in C^2(B_1)$. Iterating this process infinitely proves the result once we rescale for any $R>1$.\\

Let $\delta>0$ and $v \in \R^n$ be a unit vector. We define $\dot{u}_{\delta,v}(x, \cdot) := \frac{\dot{u}(x+\delta v, \cdot)-\dot{u}(x, \cdot) }{|\delta|^\beta}$. Note that by \eqref{eq: 5.10}, $\|\dot{u}_{\delta,v}(x, \cdot)\|_{C^\beta(B_1)} \le K$ for any choice of $\delta>0$ or $\|v\|=1$. Then, we may compute
\begin{align*}
   &\left|\frac{\dot{u}(x+\delta v, \cdot) + \dot{u}(x-\delta v, \cdot) - 2\dot{u}(x, \cdot) }{\delta^{2\beta}}\right| = \frac{1}{\delta^\beta} \left|\frac{\dot{u}(x+\delta v, \cdot)-\dot{u}(x, \cdot)}{\delta^\beta} - \frac{\dot{u}(x, \cdot)-\dot{u}(x-\delta v, \cdot)}{\delta^\beta}\right| \\
   &= \frac{\left|\dot{u}_{\delta, v}(x, \cdot) - \dot{u}_{\delta,v}(x-\delta v, \cdot) \right|}{\delta^\beta} \le K.
   \tag{5.11}\label{eq: 5.11}
\end{align*}
By \cite{stein}, this shows that $\dot{u}(x, \cdot) \in C^{2\beta}(B_1)$. We can then double regularity until $2\beta >1$ to show $\dot{u}(x, \cdot) \in C^1(B_1) $.\\

To prove \eqref{eq: 5.10}, we only need to modify the proof of \cite[Lemma 2.3]{tan} to hold for $a(t) \in A_2$ since all other steps follow independently. That is, we need to show that for $f(x,t) \in C^1_c(Q_R \cup \partial'Q_R)$:

\begin{align*}
    \int_{\partial'Q_R}|f|^2 \le \epsilon \int_{Q_R}|\nabla f|^2a(t) + \frac{C(R)}{\epsilon^{\delta}}\int_{Q_R}|f|^2a(t).
    \tag{5.12}\label{eq: 5.12}
\end{align*}
First, by calculus and Young's inequality, we may show that

\begin{align*}
    \int_{\partial'Q_R}|f|^p &= - \int_{Q_R}\partial_t|f|^p = -\int_{Q_R}p|f|^{p-1}sgn(f)\partial_t f  \\
    & \le \epsilon \int_{Q_R} |\nabla f|^p + C\epsilon^{-\frac{1}{p-1}}\int_{Q_R}|f|^p.
    \tag{5.13}\label{eq: 5.13}
\end{align*}
Next, we claim that for all $0<l<l_0$ given some $l_0>1$, we have

\begin{align*}
    \int_{Q_R}|f|^2 a(t)^{-l} \le C(l) \int_{Q_R} |\nabla f|^2 a(t).
    \tag{5.14}\label{eq: 5.14}
\end{align*}
In fact, by calculus and the Cauchy-Schwarz inequality, we compute

\begin{align*}
    f^2(x,t) = \left(\int_t^{R} \partial_t f(x,s) ds\right)^2 \le \left|\int_{t}^{R}a(s)^{-1}ds\right|\left|\int_{t}^{R}a(s)|\partial_t f|^2 ds\right|.
    \tag{5.15}\label{eq: 5.15}
\end{align*}

Since $a(t) \in A_2$, then we know that $a(t) \in A_m$ for some $1<m<2$ \cite[Lemma 5]{muckenhoupt}. By definition,

\begin{align*}
    \left(\frac{1}{R}\int_{0}^{R}a(t)\right)\left(\frac{1}{R}\int_{0}^{R}a(t)^{-\frac{m'}{m}}\right)^{\frac{m}{m'}} \le K,
    \tag{5.16}\label{eq: 5.16}
\end{align*}
where $\frac{1}{m}+\frac{1}{m'}=1$. Let $b_0 = \frac{m'}{m}$ and note that $b_0>1$ since $m<2 \implies m' >2$. Then, for $R>0$ fixed, we may assume that $\int_0^R a(t) \ge \delta_0 > 0$. Combining this fact with \eqref{eq: 5.16}, we get

\begin{align*}
    \int_0^R a^{-b_0} \le \tilde{K} \implies \int_0^R a^{-b} \le C:= R^{\frac{b_0-b}{b_0}}\tilde{K}^{\frac{b}{b_0}} \textrm{  } \forall b<b_0,
    \tag{5.17}\label{eq: 5.17}
\end{align*}
where the implication is due to H\"{o}lder's inequality. Then, since $1 < b_0$, substituting back into \eqref{eq: 5.15} yields:

\begin{align*}
    f^2(x,t) \le C\int_{t}^{R}a(s)|\nabla f|^2 ds.
    \tag{5.18}\label{eq: 5.18}
\end{align*}
Hence, multiplying by $a(t)^{-b}$ and integrating proves \eqref{eq: 5.14} since

\begin{align*}
    \int_{Q_R}a(t)^{-b}f^2 \le \tilde{C}\int_0^R a(t)^{-b}\int_{B_R}\int_0^R |\nabla f(x,s)|^2 a(s) dsdx \le C \int_{Q_R} |\nabla f|^2 a(s).
\end{align*}

Finally, using \eqref{eq: 5.13}, H\"{o}lder's and Young's inequalities, and \eqref{eq: 5.14} applied on $b=\frac{p}{2-p}$ and $b=1$ in the second last line, we show that
\begin{align*}
    & \int_{\partial' Q_R}|f|^2 = \int_{\partial' Q_R} (|f|^{\frac{2}{p}})^p
    \le \epsilon \int_{Q_R} |\nabla f^{2/p}|^p + C\epsilon^{-\frac{1}{p-1}}\int_{Q_R}|f|^2\\
    &= \epsilon (\frac{2}{p})^p \int_{Q_R} |f|^{2-p} a(t)^{-p/2}|\nabla f|^p a(t)^{p/2} + C\epsilon^{-\frac{1}{p-1}}\int_{Q_R}|f|a(t)^{-1/2}|f|a(t)^{1/2}\\
    & \le \epsilon (\frac{2}{p})^p \left(\int_{Q_R} |f|^2 a(t)^{-\frac{p}{2-p}}\right)^\frac{2-p}{2}\left(\int_{Q_R} |\nabla f|^2 a(t) \right)^{p/2} + \\
    & C\epsilon^{-\frac{1}{p-1}}\int_{Q_R} \epsilon^{1+\frac{1}{p-1}}|f|^2 a(t)^{-1} + \epsilon^{-1-\frac{1}{p-1}}|f|^2a(t)\\
    &\le \epsilon C \int_{Q_R} |\nabla f|^2 a(t) + \frac{C}{\epsilon^{1+\frac{2}{p-1}}}\int_{Q_R} |f|^2 a(t).
\end{align*}

\qed

\section{Estimates on the Harmonic Extension}

In this section, we include some new estimates for the harmonic extension $\dot{u}(x,t)$ depending on its trace $\dot{u}(x,0)=u(x)$ when $u(x)$ solves \eqref{eq: 1.5}. 

\begin{lemma}
\label{lm 6.1}
Let $\dot{u}(x,t)\in L^\infty(\R^{n+1}_+)$ solve \eqref{eq: 1.5} such that $u(x)=0$, $a(t) \in A_2$, and $\psi(-\Delta) u = u$ on $\R^n$. Then, $\dot{u}(x,t)=0$.
\end{lemma}

\noindent
\textbf{Proof: } 
Let $\tilde{a}(t)$ be the even extension of $a(t)$ so that $\tilde{a}(t)$ is an $A_2$ weight on $\R$. Then, let $\hat{a}(x,t)=\tilde{a}(t)$. $\hat{a}$ would be an $A_2$ weight on $\R^{n+1}$ by Lemma \ref{lm 6.2}. Thus, if we take the odd extension of $\dot{u}$ across $t=0$, then it will be harmonic in the sense that $\textrm{div}(\hat{a}\nabla \dot{u}) = 0$ in $\R^{n+1}$. \\

Let $Y := (Y_s)_{s \ge 0}$ be the diffusion process generated by the equation. For any open and bounded $A \subset \R^{n+1}$, let $T_A$ be the hitting time $T_A= \inf \{s \ge 0: Y_s \in A\}$ and $\tau_A$ be the exit time $\tau_A= \inf \{s \ge 0: Y_s \not \in A\}$. By the mean value property, we can express $\dot{u}(y)=\E_y[\dot{u}(Y_{\tau_A})]=\E[\dot{u}(Y_{\tau_A})| Y_0=y]$ for $y=(x,t) \in \R^{n+1}$. \\

Given $y_0 \in \R^{n+1}_+$, let $x_0$ be the projection of $y_0$ onto $t=0$. It suffices to show that $\dot{u}(y)=0$ $\forall y \in \Omega$ for any neighbourhood $\Omega = B_{R}^n(x_0) \times (0, R) \subset \R^{n+1}_+$ where $B_{R}^n(x_0) \subset \R^n$. Then, let $\Omega_{k} := B_{2^k R}^n(x_0) \times (0, 2^k R)$, $k \ge 1$. Note that $\partial \Omega_{k}$ has two parts. The first is $\Gamma_{1,k} := \{(x,t) \in \Omega_{k}: t=0\}$. Let the rest of the boundary be $\Gamma_{2,k}$. Note that $\dot{u}(y) = 0$ on $\Gamma_{1,k}$ by assumption. Hence, for any $k \ge 1$ and $y \in \Omega$, applying the mean value property yields:

\begin{align*}
    |\dot{u}(y)|=|\E_y[\dot{u}(Y_{\tau_{\Omega_k}})]| \le P_y(T_{\Gamma_{2,k}}<T_{\Gamma_{1,k}})\|\dot{u}(x,t)\|_\infty.
    \tag{6.1}\label{eq: 6.1}
\end{align*}
Since $\|\dot{u}(x,t)\|_\infty \le M <\infty$, it suffices to show that $\lim\limits_{k \to \infty}P_y(T_{\Gamma_{2,k}}<T_{\Gamma_{1,k}}) \to 0$.\\

Since $\hat{a}(t) \in A_2$ and $\textrm{div}(\hat{a}\nabla \dot{u}) = 0$, $\dot{u}$ satisfies the elliptic Harnack inequality: $\esssup_{B(y,R)} \dot{u} \le C_H \essinf_{B(y,R)} \dot{u}$ for any ball $B(y,R) \subset \R^{n+1}$ \cite[Lemma 2.3.5]{fabes}. Thus, we may apply  \cite[Lemma 3.7]{harnack}, which shows that for any $\tilde{R}>0$,

\begin{align*}
    P_y(T_{B(y_1, \tilde{R}/4)} < \tau_{B(y_0, \tilde{R})}) \ge p_0 > 0 \textrm{  }\forall y \in B(y_0, 7\tilde{R}/8),
    \tag{6.2}\label{eq: 6.2}
\end{align*}
for any $y_1 \in B(y_0, \tilde{R}/2) \subset \R^{n+1}$ and such that $p_0$ only depends on the constant $C_H$ from the Harnack inequality. Then, for any $k$, we pick $y_1^k$ as the projection of $y_0$ onto $t=-\frac{2^kR}{4}$ and let $\tilde{R}= 2^kR$. Clearly, $y_1^k \in B(y_0, \frac{2^kR}{2})$, which yields

\begin{align*}
    P_y(T_{B(y_1^k, \frac{2^kR}{4})} < \tau_{B(y_0, 2^kR)}) \ge p_0 > 0 \textrm{  } \forall y \in B\left(y_0, \frac{7*2^kR}{8}\right).
    \tag{6.3}\label{eq: 6.3}
\end{align*}

Given that the diffusion starts at $y \in B(y_0, \frac{7 \cdot 2^kR}{8})$, note that $T_{B(y_1^k, \frac{2^kR}{4})} < \tau_{B(y_0, 2^kR)}$ implies that $T_{\Gamma_{1,k}}<T_{\Gamma_{2,k}}$. Thus, for any $k$,

\begin{align*}
P_y(T_{\Gamma_{2,k}}<T_{\Gamma_{1,k}}) \le 1-p_0 \textrm{  } \forall y \in B\left(y_0, \frac{7 \cdot 2^kR}{8}\right).
\tag{6.4}\label{eq: 6.4}
\end{align*}
However, since the diffusion is continuous, then $T_{\Gamma_{2,k+1}}<T_{\Gamma_{1,k+1}}$ implies that $T_{\Gamma_{2,k}}<T_{\Gamma_{1,k}}$. Since all points on $\Gamma_{2,k}$ are in $B(y_0, \frac{7*2^{k+1}R}{8})$, it follows that:

\begin{align*}
    P(T_{\Gamma_{2,k+1}}<T_{\Gamma_{1,k+1}}| T_{\Gamma_{2,k}}<T_{\Gamma_{1,k}}) \le 1-p_0 < 1.
    \tag{6.5}\label{eq: 6.5}
\end{align*}
Then, by the strong Markov property, we obtain 
\begin{align*}
    &P_y(T_{\Gamma_{2,k}}<T_{\Gamma_{1,k}}) = P(T_{\Gamma_{2,k}}<T_{\Gamma_{1,k}}| T_{\Gamma_{2,k-1}}<T_{\Gamma_{1,k-1}})P(T_{\Gamma_{2,k-1}}<T_{\Gamma_{1,k-1}}| T_{\Gamma_{2,k-2}}<T_{\Gamma_{1,k-2}})\\
    &...P_y(T_{\Gamma_{2,1}}<T_{\Gamma_{1,1}}) \le (1-p_0)^k.
    \tag{6.6}\label{eq: 6.6}
\end{align*}
The result therefore follows after taking $k \to \infty$.

 \qed
 
 \begin{lemma}
 \label{lm 6.2}
If $a(t)$ is a non-negative $A_2$ weight on $(0, \infty)$, then it can be extended to a non-negative $A_2$ weight, $\hat{a}(x,t)$ on $\R^{n+1}$.
\end{lemma}

\noindent
\textbf{Proof:} Let $q \ge p \ge 0$. By definition $a(t) \in A_2$, $a(t)$ satisfies
\begin{align*}
    \left(\fint_p^q a(t) dt\right) \left(\fint_p^q a(t)^{-1} dt\right) \le C \iff \left(\int_p^q a(t) dt\right) \left(\int_p^q a(t)^{-1} dt\right) \le C|q-p|^2.
    \tag{6.7}\label{eq: 6.7}
\end{align*} 

Let $\tilde{a}(t)$ be the even extension of $a(t)$. That is, $\tilde{a}(-t)=\tilde{a}(t)=a(t)$ for $t \ge 0$. Then $\tilde{a}(t)$ is an $A_2$ weight on $\R$. To see why, note that the cases $q,p \ge 0$ and $q,p \le 0$  immediately reduce to \eqref{eq: 6.7}, so we only need to check the case when $\hat{p}:=-p>0$ and $q>0$. Let $M := \max \{|p|,|q|\} = \max \{\hat{p},q\} $. We verify that
\begin{align*}
     & \left(\fint_p^q \tilde{a}(t) dt\right) \left(\fint_p^q \tilde{a}(t)^{-1} dt\right)\\
     &= \frac{1}{|q-p|}\left( \int_p^0 \tilde{a}(t)dt + \int_0^q \tilde{a}(t)dt\right )\frac{1}{|q-p|}\left( \int_p^0 \tilde{a}(t)^{-1}dt + \int_0^q \tilde{a}(t)^{-1}dt\right)\\
     & \le \frac{1}{|q-p|^2} \left(\int_0^{\hat{p}} a(t)dt + \int_0^q a(t)dt \right)\left(\int_0^{\hat{p}} a(t)^{-1}dt + \int_0^q a(t)^{-1}dt \right)\\
     & \le \frac{1}{|q-p|^2} \left(C|p|^2+ 2C|M|^2 + C|q|^2 \right) \le \frac{4C|M|^2}{|M|^2}=4C.
\end{align*}
In the last line, we have expanded the product and used the non-negativity of $a(t), \frac{1}{a(t)}$ along with \eqref{eq: 6.7}. Note also that $|q-p| = q+\hat{p} \ge M$.\\

Next, we define $\hat{a}(x,t) = \tilde{a}(t)$. Let $B^n(x_0, R) \subset \R^n$ be a ball of radius $R$ centered around $x_0$, and let $t_0 \in \R$. Then, for any set $Q_R(x_0, t_0)= B^n(x_0,R) \times (t_0+\frac{R}{2},t_0-\frac{R}{2})$, the $A_2$ condition for $\hat{a}(x,t)$ is satisfied on $\R^{n+1}$ since it immediately reduces to the $A_2$ condition on $\tilde{a}(t)$. Now, let $y=(x_0,t_0) \in \R^{n+1}$. For any ball $B(y,R) \subset \R^{n+1}$, note that we may choose $\tilde{R}$ such that $B(y,R) \subset Q_{\tilde{R}}(x_0, t_0)$, and $| Q_{\tilde{R}}(x_0, t_0)| \le 2 |B(y,R)|$. Hence, we have
\begin{align*}
    & \int_{B(y,R)}\hat{a}(y)dy \le \int_{Q_{\tilde{R}}(x_0, t_0)}\hat{a}(x,t)dxdt \\
    & \le C | Q_{\tilde{R}}(x_0, t_0)|^2 \le 4C |B(y,R)|^2.
\end{align*}

\qed

 \begin{lemma}
 \label{lm 6.3}
If $\dot{u}(x,t)$ solves \eqref{eq: 1.5}, such that its trace $\dot{u}(x,0)=u(x)$ vanishes at $\infty$, then $\lim\limits_{h_n \to \infty}\|\dot{u}(x,t)\|_{Q_R(x+h_n)}=0$ where $Q_R(x) = B_R(x) \times (0,R) \subset \R^{n+1}_+$.
\end{lemma}

\noindent
\textbf{Proof:} Fix $x \in \R^n$ and let $\{h_n\}$ be a sequence in $\R^n$ such that $\|h_n\| \to \infty$. Define $\dot{u}_n(x,t) = \dot{u}(x+h_n,t)$. For each $n$, $\dot{u}_n(x,t)$ solves \eqref{eq: 1.5} with trace $u(x+h_n)$. By classical elliptic regularity theory, we know that the limit of solutions $\dot{u}^* = \lim\limits_{n \to \infty} \dot{u}_n $ is also a solution to \eqref{eq: 1.5} with trace $\lim_{n \to \infty} u(x+h_n)=0$. Therefore $\dot{u}^*$ satisfies the conditions for Lemma \ref{lm 6.1}, which concludes the proof since $ \lim\limits_{n \to \infty} \dot{u}^*_n = 0$.

\qed

\section{Conclusion}

\indent
 
In this paper, we have shown that the classical Helmholtz solutions $-\Delta u = u$ on $\R^n$ also solve many other different Helmholtz solutions. First, we classified these solutions as the bounded and vanishing fractional Helmholtz solutions for the case $0<s<1$. In the special case $n=1$, we have proved that the Helmholtz solutions given by $A\sin(x)+B\cos(x)$ are the  bounded fractional solutions using Fourier analysis. This classification extends to the case $1<s\le2$ and the polyharmonic case $s \in \N$, provided that $u \in C^\infty$. The uniqueness proofs use two techniques: the extension and energy monotonicity of \cite{Frank}, and a decomposition principle.  An open question is that this classification result can be extended for any $s \in (2,+\infty) \backslash \N$, however the decomposition technique shown in sections 3 and 4 will not suffice for reducing the case $s>2$ to the case $0<s \le 1$. It is also an open question whether or not there exist non-vanishing fractional Helmholtz solutions. 

We have also classified Helmholtz solutions for the more general equation $\psi(-\Delta)u = \psi(1)u$ when $\psi$ is completely Bernstein. The classical Helmholtz solutions solve this equation and we have proven that they are unique in the class of bounded and vanishing functions for $n \ge 2$ when we impose the regularity conditions $a(t) \in A_2$ and $a(t) \sim t^\alpha$ for some $|\alpha|<1$ and $t\gg 1$. These conditions were imposed to obtain appropriate estimates used in Lemma \ref{lm 5.2} and Lemma \ref{lm 5.4} in order to complete our uniqueness argument. Further research should be done on whether these regularity restrictions can be relaxed while preserving the uniqueness of classical Helmholtz solutions for complete Bernstein Helmholtz solutions.

\section{Appendix}

\label{sec1}

We briefly classify the solutions to the classical Helmholtz equation and note their relationship to Bessel functions for dimensions $n \ge 2$. See also \cite{Agmon} for a complete representation.

Let us first recall that the bounded Bessel function of order $\nu$, denoted by $J_\nu (r)$, satisfies the differential equation $u_{rr}+\frac{1}{r}u_r+(1-\frac{\nu^2}{r^2})u=0$.

Let $n\geq 2$. Using separation of variables and  $\Delta_{\R^n}u=u_{rr}+\frac{n-1}{r}u_r + \frac{\Delta_{\mathbb{S}^{n-1}}u}{r^2}$, it is evident that Helmholtz solutions are of the form:
\begin{center}
$u(r, \theta)=\sum\limits_{l=0}^{\infty} \sum\limits_{m=-l}^{l} c_{l,m} r^{\frac{2-n}{2}}J_{n/2+l-1}(r) \phi_{l,m}(\theta)$
    \end{center}
where     $\phi_{l,m}$ is a Laplacian spherical harmonic of order $l$ and multiplicity $m$:
\begin{center}
$\Delta_{\mathbb{S}^{n-1}}\phi_{l,m} + \mu_{l,m}\phi_l=0$, s.t. $\mu_{l,m} = l(l+n-2)$. 
\end{center}

\section*{Acknowledgement}
This research is partially supported by NSERC of Canada.

\printbibliography

\end{document}